\documentclass[10pt, a4paper, notitlepage]{article}

\usepackage{amssymb}
\usepackage{amsmath}
\usepackage{amsthm}
\usepackage{enumerate}
\usepackage{natbib}
\usepackage{cleveref}
\usepackage{color}

\newcommand{\re}{\RR}
\newcommand{\Z}{\ZZ}
\newcommand{\T}{\TT}

\renewcommand{\P}{\PP}

\newcommand{\PP}{\mathbb{P}}

\newcommand{\eps}{\varepsilon}

\newcommand{\RR}{\mathbb{R}}

\newcommand{\NN}{\mathbb{N}}

\newcommand{\BB}{\mathcal{B}}

\newcommand{\ZZ}{\mathbb{Z}}

\newcommand{\HHH}{\mathbb{H}}

\newcommand{\given}{\,|\,}
\newcommand{\TT}{\mathbb{T}}

\theoremstyle{plain}
\newtheorem{theorem}{Theorem}[section]

\newtheorem{lemma}[theorem]{Lemma}

\theoremstyle{definition}

\begin{document}

\title{\bf Gaussian process methods for one-dimensional diffusions: optimal rates and adaptation\\[.5cm]}

\author{Jan van Waaij and Harry van 
Zanten\footnote{Korteweg-de Vries Institute for Mathematics, Science Park 107, 1098 XG Amsterdam, The Netherlands.
 Email: j.vanwaaij@uva.nl, hvzanten@uva.nl.}
 \footnote{Research  funded by the
Netherlands Organization for Scientific Research (NWO).}\\[1cm]
}

\date{February 8, 2016}

\maketitle

\begin{center}
	\textit{Accepted for publication in Electronic Journal of Statistics}
\end{center}
\bigskip

\begin{abstract}
We study the  performance of nonparametric Bayes procedures
for  one-dimensional diffusions with periodic drift. We improve 
existing convergence rate results for Gaussian process (GP) priors with fixed hyper parameters.
Moreover, we exhibit several possibilities to achieve adaptation 
to smoothness. We achieve this by considering hierarchical procedures 
that  involve  either a prior on a  multiplicative scaling parameter, or 
a prior on the regularity parameter of the GP. 
\end{abstract}

\numberwithin{equation}{section}

\bigskip

\bigskip

\section{Introduction}

Various papers have recently considered nonparametric Bayes procedures 
for one-dimensional stochastic differential equations (SDEs) with periodic drift.
This is motivated among others by problems in which SDEs are used for the 
dynamic modelling of angles
in different contexts. See for instance  \cite{Hindriks} for applications in 
the modelling of neuronal rhythms and \cite{Pokern} for the use of SDEs 
in the modelling of angles in molecular dynamics.

The first paper to propose a concrete nonparametric Bayesian method in this context
 and to study its implementation was \cite{Papa}.
In \cite{PSZ} the first theoretical results were obtained for this procedure. 
These papers consider observations $(X_t: t \in [0,T])$ from  the basic SDE model
\begin{equation}\label{eq: sde}
dX_t = b(X_t)\,dt + dB_t, \qquad X_0 = 0,
\end{equation}
where $B$ is a Brownian motion, and the drift function $b$ belongs 
to the space $\dot L^2(\TT)$ of square integrable, periodic functions on 
$[0,1]$ with zero mean, i.e.\ $\int_0^1 b(x)\,dx = 0$. 
For the function $b$ of interest a GP prior is proposed with mean 
zero and precision (inverse covariance) operator 
\begin{equation}\label{eq: prior1}
\eta ((-\Delta)^{\alpha+1/2} + \kappa I),
\end{equation}
where $\Delta$ is the one-dimensional Laplacian, $I$ is the identity operator and $\eta, \kappa >0$ and 
$\alpha+1/2 \in \{2, 3, \ldots\}$ ($p=\alpha+1/2$ in \cite{PSZ}) are fixed hyperparameters. It can be 
{proved} that this  defines a valid 
prior on $\dot L^2(\TT)$, cf.\ \cite{PSZ}, Section 2.2. 

The main convergence result proved in \cite{PSZ} asserts that if in this setup 
the true drift $b_0$ generating the data has (Sobolev) regularity  $\alpha+1/2$, then 
the corresponding posterior distribution of $b$ contracts around $b_0$ at the 
rate $T^{-\alpha/(1+2\alpha)}$ as $T \to \infty$, with respect to the $L^2$-norm.
In the concluding section of \cite{PSZ} it was already conjectured that this result
is not completely sharp. More specifically, it was anticipated that the rate $T^{-\alpha/(1+2\alpha)}$
should already be attainable under the less restrictive assumption that the drift $b_0$ 
has regularity of order $\alpha$. The first main result in  the present paper confirms  that this 
is indeed the case. Since the degree of regularity of the GP with precision 
\eqref{eq: prior1} is (essentially) $\alpha$ (see e.g.\ \cite{PSZ}, Lemma 2.2.), 
this reconciles the result for this SDE model 
with the general message from the Gaussian prior literature, which says that 
to obtain optimal rates with fixed GP priors, one should match the regularities 
of the prior and the truth (see \cite{vaart2008b}).
Although lower bounds for the minimax rate appear to be unknown for the exact model we consider in this paper, 
results for closely related models suggest it is of the order 
 $T^{-\alpha/(1+2\alpha)}$ for an   $\alpha$-Sobolev smooth drift function (e.g.\ \cite{Kut}).

We are able to 
obtain {the improved result} by following a different mathematical route than in \cite{PSZ}. 
The latter paper uses more or less explicit representations of the posterior mean 
and covariance in terms of weak solutions of certain differential equations to study the asymptotic
behaviour of the posterior using techniques from PDE theory. In the present paper we  follow instead the approach of 
\cite{meulen2006}, which is essentially an adaptation to the SDE case of the general 
``testing approach'' which has by now become well known in Bayesian nonparametrics. 
These ideas, combined with results about the asymptotic behaviour of the so-called periodic 
diffusion local time from \cite{PSZ}, allow us to obtain the new, sharp result for the GP 
prior with precision \eqref{eq: prior1}.

The scope of this result is still somewhat limited, since it is a non-adaptive statement.
Indeed, it is not realistic to assume that we know the regularity of the truth  exactly and hence it 
is unlikely that we guess the correct smoothness of the prior leading to the optimal 
contraction rate. We therefore also consider several ways of obtaining 
adaptation to smoothness for this problem. A first option we explore is putting 
a prior on the multiplicative constant $\eta$ in \eqref{eq: prior1}, instead of 
taking it fixed as in \cite{Papa} and \cite{PSZ}. This  leads to 
a hierarchical, conditionally Gaussian prior on the drift $b$. Our second main result 
shows that  if the hyperprior on 
$\eta$ is appropriately chosen, then adaptation is obtained for the whole range of 
regularities between $0$ and $\alpha+1/2$. More precisely, if the degree of regularity $\beta$ of the true drift 
belongs to $(0,\alpha+1/2]$, then we attain the posterior contraction rate $T^{-\beta/(1+2\beta)}$. 

It is obviously desirable to have a large range of regularities to which we can adapt. 
At first sight, the result just discussed might suggest to let 
$\alpha$ tend to infinity with $T$. 
However, it turns out that the parameter $\alpha$ 
appears in the constant multiplying the rate of contraction.
A straightforward adaptation of the proof of the previous result  (which we will 
not carry out in this paper, since it contains no new ideas) shows that although 
taking a hyperparameter $\alpha_T \to \infty$ 
would  indeed lead to adaptation over the growing interval $(0,\alpha_T+1/2]$, 
the rate would deteriorate by a factor $(\alpha_T)^c$ for some constant $c > 0$. 

The preceding observations  indicate that in  order to obtain adaptation to the full range of 
possible regularities for the drift, using a prior on the multiplicative scale parameter $\eta$
is perhaps not the best option. Therefore we  also consider another possibility, 
namely putting a prior on the hyperparameter $\alpha$ that controls the regularity of the 
prior directly. We prove that this is, from the theoretical perspective at least, indeed
preferable. We can obtain the optimal 
contraction rate for any regularity of the truth, without suffering a penalty in the rate.

In this paper we focus on deriving theoretical results. 
We do not consider the related numerical issues, since this requires a completely different 
analysis, but these are clearly of interest as well. 
For instance, it is quite conceivable 
that the last option we consider, putting a prior on $\alpha$, is numerically quite demanding, 
more so than putting a prior on $\eta$. Therefore in practice it might actually be worthwhile 
to accept non-optimal statistical rates or only a limited range of adaptation,  in order to gain speed 
on the numerical side. 
The paper \cite{Moritz} considers a related but different computational strategy, 
which combines a prior on the multiplicative constant with a random truncation 
of the series that defines the Gaussian prior. 
Related to this is the work of  \cite{agapiou2014}, who study similar approaches
in different statistical settings. It would be of interest to understand 
the theoretical performance of such computationally attractive methods better. This is outside the scope 
of the present paper however and remains to be dealt with in forthcoming work.

The paper is organised as follows. In the next section we describe the diffusion model 
and the priors that we consider  in detail. 
In Section \ref{sec: main} we present and discuss the main results described briefly in the introduction. 
Some auxiliary result that we use in the proofs are prepared in Section \ref{sec: aux}. 
The proofs themselves are given in Sections \ref{sec: proof1}--\ref{sec: proof3}.

\section{Model and prior}
\label{sec: model}

As explained in the introduction we consider the 1-periodic diffusion model given by 
\eqref{eq: sde}, where $B$ is a standard Brownian motion 
and $b:\re\to\re$ is a measurable function that is 1-periodic, square integrable and mean zero on $[0,1]$. 
 The space of all such functions is denoted by $\dot L^2(\T).$ We endow this space with the 
 usual $L^2$-norm defined by $\|b\|^2_2=\int_0^1b^2(x)\,dx$.
We note that for any $b \in \dot L^2(\TT)$, the SDE \eqref{eq: sde} admits 
a unique weak solution. (For the sake of completeness 
we have added a proof in the appendix, see Lemma \ref{lemma:uniquesolutiontosde}.)

For every $T > 0$ the solution $X^T = (X_t: t \in [0,T])$ of the SDE 
induces a law $P^{b}=P^{b,T}$ on the space $C[0,T]$ of continuous functions on $[0,T]$. 
For fixed $T>0,$ $P^{b_1}$ and $P^{b_2}$ are equivalent for all $b_1,b_2\in\dot L^2(\T)$ (see Lemma \ref{lemma:equivalentmeasures}). 
The Radon-Nikodym derivative $p_b$ of $P^{b,T}$ relative to the Wiener measure ($P^{0,T}$) satisfies 
\begin{equation}\label{eq: lik}
p_b(X^T) = \frac{d P^{b,T}}{d P^{0,T}}(X^T)=\exp\left(-\frac{1}{2}\int_0^T b^2(X_t)dt+\int_0^Tb(X_t)dX_t\right)
\end{equation}
almost surely.

To make Bayesian inference about the drift function we consider a Gaussian process (GP) prior on the space 
of drift functions $\dot L^2(\TT)$. We are interested in the GP with 
mean zero and precision operator \eqref{eq: prior1}. 
As shown in Section 2.2 of \cite{PSZ}, the GP $W$ with this mean and 
covariance can be written as 
\[
W = \frac1{\sqrt\eta}\sum_{k=1}^\infty \sqrt\lambda_k \phi_k Z_k,  
\]
where the $Z_k$ are independent standard normal variables, the $\phi_k$ are the orthonormal eigenfunctions
of the Laplacian, given by 
\begin{align*}
\phi_{2k}(x) & = \sqrt 2 \cos(2\pi k x),\\
\phi_{2k-1}(x) & = \sqrt 2 \sin(2\pi k x),
\end{align*}
for $k \in \NN$, and 
\begin{equation}\label{eq: l1}
\lambda_k = \Big(\Big(4\pi^2\Big\lceil \frac k2 \Big\rceil^2\Big)^{\alpha+1/2}+ \kappa\Big)^{-1}.
\end{equation}

The results we derive in this paper actually do not depend crucially on the exact 
form of the eigenfunctions and eigenvalues $\phi_k$ and $\lambda_k$. 
The $\phi_k$ can in fact be {any} orthonormal basis of $\dot L^2(\TT)$
(provided the smoothness spaces defined ahead are changed accordingly). 
Moreover, the specific value of the hyperparameter $\kappa$ in \eqref{eq: l1}
is irrelevant for our results.
For the $\lambda_k$ we only need that there exist constants $c, C > 0$ and $\alpha > 0$ 
such that 
\begin{equation}\label{eq: l2}
ck^{-1/2-\alpha} \le \sqrt\lambda_k \le Ck^{-1/2-\alpha}. 
\end{equation}
Note that the $\lambda_k$'s in \eqref{eq: l1} satisfy these bounds. 
For notational convenience we will work with $\sqrt\lambda_k = k^{-1/2-\alpha}$ 
throughout the paper, but all results hold if this exact choice is replaced by 
$\lambda_k$'s satisfying \eqref{eq: l2}.

Introducing the notation $L = 1/\sqrt\eta$ for the scaling constant, the priors 
for the drift function $b$ that we consider take the general form 
\begin{equation}\label{eq: p}
b \sim L \sum_{k=1}^\infty k^{-1/2- \alpha}\phi_k Z_k, 
\end{equation}
where the $Z_k$ are independent standard Gaussian variables and  $(\phi_k)$ is 
an arbitrary, fixed orthonormal basis of $\dot L^2(\TT)$. We will consider  various 
setups in which 
the scale $L$ is either a constant 
or a random factor, and also the regularity parameter $\alpha$ will be either a fixed constant
or random. 

The regularity of the true drift function $b_0$ that generates the data
will be measured in Sobolev sense relative to the basis $(\phi_k)$. 
For $\beta > 0$ we define 
\[
\dot H^\beta(\TT) = \Big\{f = \sum_{k=1}^\infty f_k\phi_k \in \dot L^2(\TT): \sum_{k=1}^\infty f_k^2 k^{2\beta} < \infty\Big\}. 
\]
Note that  in the case that the $\phi_k$ are the eigenfunctions of the Laplacian 
given above, this is the usual $L^2$-Sobolev regularity.

\section{Main results}
\label{sec: main}

In this section we present the main rate of contraction results for 
the posteriors corresponding to the various priors of the form \eqref{eq: p}, with different choices
for the hyperparameters $L$ and $\alpha$. The proofs of the results are given in \Cref{sec: proof1,sec: proof2,sec: proof3}.

For simplicity the prior on $b$ will always be denoted by $\Pi$, but it will be clearly described in each
case. For every time horizon $T > 0$, the corresponding posterior distribution will be 
denoted by $\Pi(\cdot \given X^T)$.
So for  a Borel set $A \subset \dot L^2(\TT)$, 
\[
\Pi(b \in A\given X^T) = \frac{\int_A p_b(X^T)\,\Pi(db)}{\int p_b(X^T)\,\Pi(db)},
\]
where the likelihood is given by \eqref{eq: lik}.
The following lemma asserts that the posterior is well defined under the minimal condition
that the prior $\Pi$ is a probability measure on the Borel sets of $\dot L^2(\TT)$. 
The proof is deferred to \Cref{proof:posteriorwelldefined}.

\bigskip

\begin{lemma}\label{lem:posterioriswelldefined} Suppose that $\Pi$ is  Borel probability measure on $\dot L^2(\TT)$.  
Then for every $b_0 \in \dot L^2(\TT)$ it $P^{b_0}$-a.s.\ holds that 
\begin{enumerate}
\item[(i)]
the random map $b \mapsto p_b(X^T)$ admits a version that is Borel measurable on $\dot L^2(\TT)$, 
\item[(ii)]
for the denominator we have $0<\int p_b(X^T)\,\Pi(db) <\infty$.
\end{enumerate}
\end{lemma}

\bigskip

 As usual we say that 
the posterior {\em contracts around $b_0$ at the rate $\eps_T$ as $T \to \infty$}
if for all $M_T \to \infty$, 
\[
\Pi\Big(b: \|b - b_0\|_2 \ge M_T\eps_T \given X^T\Big) \overset{P^{b_0}}{\longrightarrow} 0 
\]
as $T \to \infty$. Here the convergence is in probability under the law $P^{b_0}$
corresponding to the true drift function $b_0$.

\subsection{Fixed hyperparameters}

Our first main result deals with the case that the scaling parameter  and the 
regularity parameter of the GP are fixed, positive constants. Specifically, 
we fix $L > 0$ and $\beta > 0$ and define the prior $\Pi$ on the drift function 
structurally  as
\begin{equation}\label{eq: pp1}
b \sim L \sum_{k=1}^\infty k^{-1/2- \beta}\phi_k Z_k, 
\end{equation}
where the $Z_k$ are independent standard Gaussian variables and  $(\phi_k)$ is 
the chosen orthonormal basis of $\dot L^2(\TT)$. 
Note that the expected squared $L^2$-norm of $b$ under this prior is $L^2\sum k^{-1-2\beta} < \infty$, 
hence by Lemma \ref{lem:posterioriswelldefined} the posterior is well defined.

\bigskip

\begin{theorem}\label{thm: 1}
Let the prior be given by \eqref{eq: pp1}, with $\beta, L > 0$ fixed. 
If $b_0 \in \dot H^\beta(\TT)$ for $\beta > 0$, then the posterior contracts around $b_0$
at the rate $\eps_T = T^{-\beta/(1+2\beta)}$. 
\end{theorem}

\bigskip

As noted in the introduction, this theorem improves  Theorem 5.2 of \cite{PSZ}. The latter corresponds to 
the case that the $\phi_k$ are the eigenfunctions of the Laplacian and $\beta+1/2 \in \{2, 3, \ldots\}$.  
In \cite{PSZ} the obtained rate for this prior is also (essentially) $T^{-\beta/(1+2\beta)}$, but 
this is obtained under the stronger condition  that $b_0$ belongs to $\dot H^{\beta+1/2}(\TT)$. 
Additionally, the new result is valid for all $\beta > 0.$

\subsection{Prior on the scale}

The fact that we get the optimal rate $T^{-\beta/(1+2\beta)}$ in Theorem \ref{thm: 1} 
strongly depends on the fact that the degree of smoothness $\beta$ of the true 
drift $b_0$ matches the choice of the regularity parameter of the prior. 
Although strictly speaking it has not been established for the SDE setting of this paper, 
results from the GP prior literature for analogous settings indicate that 
if these regularities are not matched exactly, then sub-optimal rates will 
 be obtained (see for instance \cite{vaart2008b} and \cite{ismael}). 
We would obviously prefer a method that does not depend on knowledge of the 
true regularity $\beta$ of the truth and that {\em adapts} to this degree of smoothness
automatically. 

In this section we consider a first method to achieve this. This involves putting 
a prior distribution on the scaling parameter $L$ instead of taking it fixed. 
We employ a  hierarchical prior $\Pi$ on $b$ that can be described as follows:
\begin{align}
\label{eq: pp21} L & \sim \frac{E^{1/2+\alpha}}{\sqrt T},\\
\label{eq: pp22} b \given L & \sim  L \sum_{k=1}^\infty k^{-1/2- \alpha}\phi_k Z_k.
\end{align}
Here $\alpha > 0$ is a fixed hyperparameter, which should be thought of 
as describing the ``baseline smoothness'' of the prior. The $Z_k$ and $\phi_k$ are as before
and $E$ is a standard exponential, independent of the $Z_k$. Note that we could 
equivalently describe the prior on $L$ as a Weibull distribution with scale parameter $1/\sqrt T$
and shape parameter $2/(1+2\alpha)$.  Lemma \ref{lem:posterioriswelldefined} ensures again 
that the posterior is well defined, since by conditioning we see that the expected squared $L^2$-norm
of $b$ is now given by $c_L\sum k^{-1-2\alpha}$, where $c_L$ is the second moment of $L$
under the prior, which is finite.

The specific choice of the prior for $L$ is convenient, but the proof of the following theorem 
shows that it can actually be slightly generalised. It is for instance enough that the random 
variable $E$ in \eqref{eq: pp21} has a density that satisfies exponential lower and upper bounds
in the tail. Our proof breaks down however if we deviate too much from the choice above. 
For instance, without the dependence on $T$ we would only be able to derive sub-optimal rates. 
We stress that this does not mean that other priors cannot lead to optimal rates, 
only that such results cannot be obtained using our technical approach. An alternative 
route, for instance via empirical Bayes as in \cite{Bartek}, might lead to less restrictive 
assumptions on the hyperprior for $L$. This will require a completely different analysis however.

\bigskip

\begin{theorem}\label{thm: 2}
Let the prior be given by \eqref{eq: pp21}--\eqref{eq: pp22}, 
 with  $\alpha > 0$ fixed. If $b_0 \in \dot H^\beta(\TT)$ for $\beta \in (0, \alpha + 1/2]$, then the posterior contracts around $b_0$
at the rate $\eps_T = T^{-\beta/(1+2\beta)}$. 
\end{theorem}

\bigskip

So indeed with a prior on the multiplicative scale we can achieve adaptation for a range of 
smoothness levels $\beta$. 
Note however that the 
range is limited by the baseline smoothness $\alpha$ of the prior. Putting a prior on the scale $L$
does allow to adapt to truths that are arbitrarily rougher than the prior, but if the degree of smoothness of 
the truth is larger than $\alpha + 1/2$, the procedure does not achieve optimal rates. 
This phenomenon has been observed in the literature in 
different statistical settings as well. See for instance \cite{SVZ} for similar 
results in the white noise model.

\subsection{Prior on the GP regularity}

To circumvent the potential problems described in the preceding section, we 
consider an alternative method for achieving adaptation to all smoothness levels. 
Instead of taking a fixed baseline prior smoothness and putting a prior on the scale, 
we put a prior on the GP smoothness itself.
Specifically, we use a prior on $\alpha$
that is truncated to the growing interval $(0, \alpha_T]$ and that has a density proportional to $x \mapsto \exp(-T^{1/(1+2x)})$ on that interval. For convenience we take $\alpha_T=\log T,$ but other choices are possible as well. 
We define the probability density $\lambda_T$,  with support $[0, \log T]$,  
 by 
\[
\lambda_T(x) = C_T^{-1}{e^{-T^{1/(1+2x)}}}, \qquad x \in [0,\log T],
\]
where $C_T$ is the normalising constant. 
The full prior $\Pi$ on $b$ that we employ is now described as follows:
\begin{align}
\label{eq: pp31} \alpha & \sim \lambda_T,\\
\label{eq: pp32} b \given \alpha & \sim \sum_{k=1}^\infty k^{-1/2- \alpha}\phi_k Z_k,
\end{align}
where the $\phi_k$ and $Z_k$ are again as before. 
Note that for this prior we have that for every $\alpha > 0$, the conditional 
prior probability that $\|b\|_2 < \infty$ given $\alpha$ equals $1$, hence the unconditional 
prior probability that the norm is finite is $1$ as well. Lemma \ref{lem:posterioriswelldefined}
thus implies the posterior is well defined again and we can formulate the 
following result.

\bigskip

\begin{theorem}\label{thm: 4}
Let the prior be given by \eqref{eq: pp31}--\eqref{eq: pp32}. If $b_0 \in \dot H^\beta(\TT)$ for $\beta> 0$, then the posterior contracts around $b_0$
at the rate $\eps_T =  T^{-\beta/(1+2\beta)}$. 
\end{theorem}

\bigskip

So by placing a prior on $\alpha$ we obtain adaptation to all smoothness levels, 
without paying for it in the rate. A similar result has recently been obtained
in the setting of the white noise model in \cite{Bartek}. We note however 
that the results in the latter paper rely on rather explicit computations specific 
for that model. The results we present here for the SDE model are derived in a completely different way, 
by using the testing approach proposed in \cite{meulen2006}. 
We note that the rates we obtain are slightly better than those in  \cite{Bartek}, 
in the sense  that we don't obtain additional slowly varying factors. 
We expect that similar results can be obtained for white noise model and other related models
by adapting our proofs.

A downside of  our approach 
 is that we can only prove the desired result for somewhat contrived 
hyperpriors on $\alpha$ such as $\lambda_T$, which may appear unnatural at first sight.
The result is however in accordance with similar findings for other statistical models 
obtained for instance in \cite{lember2007} and \cite{ghosallembervaart2008}.
Our prior on $\alpha$ has a density proportional (on $(0,\alpha_T]$) to 
$\exp(-T\eps^2_{\alpha, T})$, where $\eps_{\alpha, T}$ is the rate we would 
get when using the unconditional Gaussian prior on the right of \eqref{eq: pp32}.
Hence our theorem is in accordance with the results in the cited papers, 
which state that in some generality, such a choice of hyper prior leads to 
rate-adaptive procedures. Other priors on $\alpha$ may lead to adaptation as well,
including potentially priors that do not depend on the sample length $T$.  
But to prove such results, different mathematical techniques seem to be required.

The main point we want to make here 
however, and that is supported by the theorems we present,  is 
that if the goal is to achieve adaptation to 
an unrestricted range of smoothness levels, then, from the theoretical point 
of view at least, putting a prior on a smoothness hyperparameter is preferable 
to fixing the baseline smoothness of the prior and putting a prior on a multiplicative 
scaling parameter.

\section{Auxiliary results}
\label{sec: aux}

Here we prepare a number of results that will be used in the proofs of Theorems 
\ref{thm: 1}--\ref{thm: 4}.

\subsection{General contraction rate result}

In this section we first present a  contraction rate result 
for general posteriors in the setting of one-dimensional SDEs 
with periodic drift, as described in Section \ref{sec: model}. 
This theorem is a consequence of the general result
of \cite{meulen2006}, in combination with a result on the 
periodic local time of the solution to \eqref{eq: sde}. 
The result is in the spirit of the corresponding i.i.d.\
result of \cite{GGV} and gives conditions for having a certain 
rate of contraction in terms of the prior mass around the truth, 
and the complexity of the essential support of the prior. 
In the sections ahead we apply it to the priors considered in Section \ref{sec: main}. 

The prior in the general theorem may depend on the 
time horizon $T > 0$ and is denoted by $\Pi_T$. For a metric space $(A, d)$
and $\eps > 0$, we denote by $N(\eps, A, d)$ the minimal number of balls of 
$d$-radius $\eps$ needed to cover the set $A$. Recall that we 
say that the posterior contracts around $b_0$ at the rate $\eps_T$ as $T \to \infty$
if for all $M_T \to \infty$, 
\[
\Pi_T\Big(b: \|b - b_0\|_2 \ge M_T\eps_T \given X^T\Big) \overset{P^{b_0}}{\longrightarrow} 0 
\]
as $T \to \infty$. 

\bigskip

\begin{theorem}\label{thm: general}
Let $\eps_T \to 0$ be positive numbers such that $T\eps_T^2 \to \infty$.  
Suppose that for some $C_1 > 0$, 
\begin{equation}\label{eq: pm}
\Pi_T(b: \|b - b_0\|_2 \le  \eps_T) \ge e^{-C_1 T \eps^2_T}.
\end{equation}
Moreover, assume that for any $C_2 > 0$, there exist measurable subsets $\BB_T \subset \dot L^2(\TT)$ 
and a $C_3 > 0$ such that
\begin{align*}
\Pi_T(\BB_T) & \ge 1-e^{-C_2T\eps^2_T},\\
\log N(\eps_T, \BB_T, \|\cdot\|_2) & \le C_3T\eps^2_T. 
\end{align*} 
Then the posterior contracts around $b_0$ at the rate $\eps_T$ as $T \to \infty$. 
\end{theorem}

\bigskip

\begin{proof}
The result follows from Theorem 2.1 and Lemma 2.2 of \cite{meulen2006}, provided that 
we show, in accordance with Assumption 2.1 of the latter paper, 
that the random distance whose square is given by 
\[
\frac1{T}\int_0^T (b(X_t)- b_0(X_t))^2\,dt,
\]
is with $P^{b_0}$-probability tending to $1$ equivalent to the $L^2$-norm $\|b-b_0\|_2$. 
But this easily follows from the asymptotic properties of the so-called 
periodic local time $(L^\circ_t(x; X): t \ge 0, x \in [0,1])$ of the process $X$ derived in 
\cite{PSZ}.

Indeed, by the occupation times formula for the periodic local time the 
integral in the preceding display equals 
\[
\int_0^1 (b(x)-b_0(x))^2 \frac1TL^\circ_T(x; X)\,dx,  
\]
see Section 2.1 of \cite{PSZ}. By the uniform law of large numbers 
given in Theorem 4.1.(i) of \cite{PSZ}, the random function $L^\circ_T/T$
converges uniformly to the invariant density $\rho$ on $[0,1]$ with $P^{b_0}$-probability $1$, 
which is given by 
\[
\rho(x) = Ce^{2\int_0^x b_0(y)\,dy}, \qquad x \in [0,1], 
\]
where $C > 0$ is the normalising constant. Since $\rho$ is bounded 
away from $0$ and $\infty$ on $[0,1]$, this shows that for every $\gamma > 0$, 
there exist constants $C_1, C_2 > 0$ such that with 
$P^{b_0}$-probability at least $1-\gamma$, and for all $b \in L^2(\TT)$,
\[
C_1\|b-b_0\|^2_2 \le \frac1{T}\int_0^T (b(X_t)- b_0(X_t))^2\,dt \le C_2 \|b-b_0\|^2_2.
\]
This is the desired equivalence of norms.
\end{proof}

\subsection{Small ball probabilities}
\label{sec: sb}

In this section we prepare a result that allows us to verify the 
prior mass condition \eqref{eq: pm} of Theorem \ref{thm: general} for the various priors
in Section \ref{sec: main}. For $\alpha, L > 0$ we define the GP
\begin{equation}\label{eq: w}
W^{\alpha, L} = L \sum_{k=1}^\infty k^{-1/2- \alpha}\phi_k Z_k, 
\end{equation}
where the $Z_k$ are independent standard Gaussian variables and  $(\phi_k)$ is 
an arbitrary orthonormal basis of $\dot L^2(\TT)$.

\bigskip

\begin{lemma}\label{lem:smallballbalphaL}
There exists a positive, continuous function $f$ on $(0,\infty)$  and constants $c_0, c_1 > 0$ 
such that $c_0 \alpha \le f(\alpha) \le c_1 \alpha$ for $\alpha$ large enough  and 
\[
-\log\PP(\|W^{\alpha,L}\|_2<\eps)\le f(\alpha)\Big(\frac{L}{\eps}\Big)^{1/\alpha},
\] 
for all $\alpha>0$ and for ${\eps}/{L} > 0 $ small enough. 
 \end{lemma}

\bigskip 
 
  \begin{proof}
 Note that $\P(\|W^{\alpha,L}\|_2<\eps)=\P(\|W^{\alpha,1}\|_2<\eps/L)$, 
 so the case $L=1$ implies the general case. 
Since $(\phi_k)$ is an orthonormal basis, 
 $\P(\|W^{\alpha,1}\|_2<\eps)=\P\left(\sum_{k=1}^\infty k^{-2\alpha-1}Z_k^2<\eps^2\right)$. 
 The result  then follows from Corollary 4.3 of \cite{dunker1998} and  
 straightforward algebra. 
\end{proof}

\bigskip

Next we consider the reproducing kernel Hilbert space (RKHS) $\HHH^{\alpha, L}$ associated to the GP $W^{\alpha, L}$. 
It follows from the series representation \eqref{eq: w} that 
$\HHH^{\alpha, L} = \dot H^{1/2+\alpha}(\TT)$, 
and that the associated RKHS norm of an element $h \in \HHH^{\alpha, L}$ satisfies 
$L\|h\|_{\HHH^{\alpha, L}} = \|h\|_{2, 1/2+\alpha}$,
where for $\beta > 0$, the Sobolev norm $\|h\|_{2, \beta}$ 
of a function $h = \sum h_k\phi_k$ 
is defined by  
\[
\|h\|^2_{2, \beta} = \sum_{k=1}^\infty h_k^2k^{2\beta}.
\]
For these facts and more general background on RKHS's of GP's with a view 
towards Bayesian nonparametrics, see \cite{vaart2008}.

\bigskip

\begin{lemma}\label{lem:L2normtoRKHSnorminfiniteseries} 
Suppose that $b_0 \in \dot H^\beta(\TT)$ for  $\beta\le\alpha+1/2$. 
Then for $\eps > 0$ small enough, 
\[
\inf_{h\in\HHH^{\alpha,L}:\|h-b_0\|_2\le\eps}\|h\|_{\HHH^{\alpha,L}}^2\le\|b_0\|_{2,\beta}^2\frac{1}{L^2}\eps^{\tfrac{2\beta-2\alpha-1}{\beta}}.
\]
\end{lemma}

\bigskip

\begin{proof}
Consider the expansion $b_0 = \sum_{k=1}^\infty b_k\phi_k$ and define 
$h = \sum_{k \le I} b_k\phi_k$, where $I$ will be determined  below. 
We have that $h \in \HHH^{\alpha, L}$, and from the smoothness condition 
on $b_0$ it follows that 
\begin{align*}
\|h-b_0\|_2^2 & = \sum_{k > I} b^2_k \le I^{-2\beta}\sum_{k >I} b^2_kk^{2\beta}. 
\end{align*}
Since $b_0 \in \dot H^\beta(\TT)$ the sum on the right vanishes for $I \to \infty$, 
hence $\|h-b_0\|_2^2 \le I^{-2\beta}$ for $I$ large enough.
Setting $I = \eps^{-1/\beta}$ we obtain that, for $\eps$ small enough,  
the infimum in the statement of the lemma is bounded by 
\begin{align*}
\frac 1{L^2}\sum_{k\le I} b_k^2k^{1+2\alpha} =  \frac 1{L^2}\sum_{k\le I} b_k^2k^{2\beta}k^{1+2\alpha-2\beta}
 \le \frac 1{L^2} \|b_0\|^2_{2,\beta} I^{1+2\alpha-2\beta} ,  
\end{align*}
since $\beta\le\alpha+1/2$. The proof is completed by recalling the choice of $I$.
\end{proof}

\bigskip

Lemmas \ref{lem:smallballbalphaL} and \ref{lem:L2normtoRKHSnorminfiniteseries} 
together give a non-centered small ball probability bound for the GP $W^{\alpha, L}$. 
This will be used to verify the prior mass condition \eqref{eq: pm} of Theorem \ref{thm: general}
for the various priors.

\bigskip

\begin{lemma}\label{lem: sb}
Suppose that $\alpha > 0$ and $b_0 \in \dot H^\beta(\TT)$ for  $\beta\le\alpha+1/2$. 
There exist a constant $C > 0$, depending only on  $b_0$, such that 
\[
\PP(\|W^{\alpha,L} - b_0\|_2<\eps)\ge \exp\left(- C\left(f(\alpha)\left(\frac{L}{\eps}\right)^{1/\alpha}
+\frac{1}{L^2}\eps^{\tfrac{2\beta-2\alpha-1}{\beta}}\right)\right).
\] 
for  ${\eps}/{L} > 0 $ small enough. 
\end{lemma}

\bigskip

\begin{proof}
This follows directly from Lemmas \ref{lem:smallballbalphaL} and \ref{lem:L2normtoRKHSnorminfiniteseries}
using, for instance,  Lemma 5.3 of \cite{vaart2008}.
\end{proof}

\section{Proof of Theorem \ref{thm: 1}}
\label{sec: proof1}

In this case the prior $\Pi$ is the law of GP $W^{\beta, L}$. Applying Lemma 
\ref{lem: sb} with $\alpha =\beta$ we obtain, for $b_0 \in \dot H^\beta(\TT)$, the bound 
\[
\Pi(b: \|b - b_0\|_2 \le \eps) \ge e^{-C \eps^{-1/\beta}},
\]
for a constant $C > 0$ and $\eps > 0$ small enough. It follows that the prior mass 
condition \eqref{eq: pm} of Theorem \ref{thm: general} is satisfied for $\eps_T$ a constant times 
$T^{-\beta/(1+2\beta)}$. By the general result for Gaussian priors given by Theorem 2.1
of \cite{vaart2008b}, the other assumptions of Theorem \ref{thm: general} 
are then automatically satisfied as well. Hence, the desired result follows from an application of that theorem.

\section{Proof of Theorem \ref{thm: 2}}
\label{sec: proof2}

We will again verify the conditions of Theorem \ref{thm: general}. We note that in this 
case, the conditional distribution of $b$ under the prior, given the value of $L$, 
is the law of $W^{\alpha, L}$.  

\subsection{Prior mass condition}

Denoting the prior density of $L$ by $g$, and assuming again that $b_0 \in \dot H^\beta(\TT)$, we have, by Lemma \ref{lem: sb}, 
that there exists a constant $C > 0$ such that for $\eps$ small enough,  
\begin{align*}
\Pi(b: \|b - b_0\|_2 \le \eps) & = \int \PP(\|W^{\alpha, L} - b_0\|_2 \le \eps)g(L)\,dL \\
& \ge \int_{\eps^{(\beta-\alpha)/\beta}}^{2 \eps^{(\beta-\alpha)/\beta}} e^{-C((L/\eps)^{1/\alpha} + \eps^{(2\beta-2\alpha -1)/\beta}/L^2)}g(L)\,dL.
\end{align*}
 On the range of integration the exponential in the integrand is bounded from below by 
$e^{-C'\eps^{-1/\beta}}$ for some $C' > 0$. Moreover, the assumptions on the prior on $L$ imply that 
for $\eps$ a multiple of $T^{-\beta/(1+2\beta)}$, 
\[
\int_{\eps^{(\beta-\alpha)/\beta}}^{2\eps^{(\beta-\alpha)/\beta}} g(L)\,dL = 
\PP(cT^{1/(1+2\beta)} < E < 2cT^{1/(1+2\beta)}) \ge  e^{-3cT^{1/(1+2\beta)}}
\]
for a constant $c > 0$ and $T$ large enough. It follows that there exist constants
$c_1, c_2 > 0$ such that for $\eps_T = c_1T^{-\beta/(1+2\beta)}$, 
\[
\Pi(b: \|b - b_0\|_2 \le \eps_T) \ge  e^{-c_2T\eps_T^2},
\]
which covers the first condition of Theorem \ref{thm: general}.

\subsection{Sieves}

Recall from Section \ref{sec: sb} that the RKHS unit ball $\HHH^{\alpha, L}_1$ of $W^{\alpha, L}$ is the 
ball $\dot H_L^{\alpha+1/2}(\TT)$ of radius $L$ in the Sobolev space $\dot H^{\alpha + 1/2}(\TT)$
of regularity $\alpha +1/2$. 
This motivates the definition of sieves $\BB_T$ of the form
\[
\BB_T = R \dot H_1^{\alpha+1/2}(\TT) + \eps_T\dot L_1^2(\TT),
\]
where $R$ will be determined below and $\dot L_1^2(\TT)$ is the unit ball in $\dot L^2(\TT)$.

\subsubsection{Remaining mass condition}

By conditioning we have, for any $L_0 > 0$, 
\begin{equation}\label{eq: rmts}
\begin{split}
\Pi(b \not \in \BB_T) & = \int \PP(W^{\alpha, L} \not \in \BB_T)g(L)\,dL \\
& \le \int_0^{L_0} \PP(W^{\alpha, L} \not \in \BB_T)g(L)\,dL + 
\int_{L_0}^\infty g(L)\,dL.
\end{split}
\end{equation}
The second term on the right is bounded by $\exp(-(L^2_0T)^{1/(1+2\alpha)})$, by the assumptions 
on the prior on $L$. For $L_0$ a large enough multiple of $T^{(\alpha-\beta)/(1+2\beta)}$ 
this is bounded by $e^{-DT^{1/(1+2\beta)}}$, for a  given constant $D > 0$.

As for the first term, note that the probability in the integrand is increasing in $L$.
Since $\BB_T = (R/L_0)\HHH^{\alpha, L_0}_1+ \eps_T\dot L_1^2(\TT)$, the Borell-Sudakov inequality (see
\cite{vaart2008}, Theorem 5.1)
implies that 
\[
\PP(W^{\alpha, L_0} \not \in \BB_T) \le 1-\Phi(\Phi^{-1}(\PP(\|W^{\alpha, L_0}\|_2 \le \eps_T)) + R/L_0).
\]
By Lemma \ref{lem:smallballbalphaL}, the probability on the right is bounded from below 
by $\exp(-C(L_0/\eps_T)^{1/\alpha})$ for some $C > 0$.
Furthermore, since for $y\in(0,0.5),\Phi^{-1}(y)\ge-\sqrt{\tfrac{5}{2}\log(1/y)}$ and for $x\ge1,1-\Phi(x)\le \exp(-x^2/2),$ we have 
\[
\PP(W^{\alpha, L_0} \not \in \BB_T) \le \exp\left(-\frac12\left(\frac{R}{L_0}-\sqrt{C'\Big(\frac{L_0}{\eps_T}\Big)^{1/\alpha}}\right)^2\right), 
\]
for some $C' > 0$. The choices of $L_0$ and $\eps_T$ imply that 
if $R$ is chosen to be a large multiple of $T^{(1/2+\alpha - \beta)/(1+2\beta)}$, then the 
first term on the right of \eqref{eq: rmts} is bounded by $e^{-DT^{1/(1+2\beta)}}$ as well.

\subsubsection{Entropy}

It remains to verify that $\BB_T$ satisfies the entropy condition of Theorem \ref{thm: general}. 
By the known entropy bound for Sobolev balls  we have 
\[
\log N(\eps, R \dot H^{\alpha + 1/2}_1, \|\cdot\|_2) \le C \Big(\frac R\eps\Big)^{2/(1+2\alpha)}
\]
for some $C > 0$. Recalling the definitions of $\BB_T$, $\eps_T$ and $R$, it follows that 
\[
\log N(2\eps_T, \BB_T, \|\cdot\|_2) \le C \Big(\frac R{\eps_T}\Big)^{2/(1+2\alpha)}
\le C' T^{1/(1+2\beta)}
\]
for some $C' > 0$. This concludes the proof of the theorem.

\section{Proof of Theorem \ref{thm: 4}}
\label{sec: proof3}

Note that in this case the conditional prior law of $b$, given $\alpha$, is the law 
of the GP $W^{\alpha, 1}$. 

\subsection{Prior mass condition}

By Lemma \ref{lem: sb}, there exist a constant $C > 0$ such that for $\eps$ small enough,  $\delta > 0$
and $b_0 \in \dot H^\beta(\TT)$, 
\begin{align*}
\Pi(b: \|b - b_0\|_2 \le \eps) & \ge \int_\beta^{\beta+\delta} \PP(\|W^{\alpha, 1} - b_0\|_2 \le \eps)\lambda_T(\alpha)\,d\alpha \\
& \ge \int_\beta^{\beta+\delta} e^{-C((1/\eps)^{1/\alpha} + \eps^{(2\beta-2\alpha -1)/\beta})}\lambda_T(\alpha)\,d\alpha.
\end{align*}
On the range of integration the exponential in the integrand is bounded from below by 
$\exp(-C'\eps^{-(1+2\delta)/\beta})$ for some $C' > 0$. Since $\lambda_T$ is increasing, we get
\[
\Pi(b: \|b - b_0\|_2 \le \eps) \ge \delta C_T^{-1}e^{-T^{/(1+2\beta)}} e^{-C'\eps^{-(1+2\delta)/\beta}}.
\]
Since $C_T\le\log T$ and by choosing $\delta$ to be a multiple of $1/\log T$, it follows that, for $\eps_T$ a multiple 
of $T^{-\beta/(1+2\beta)},$ condition \eqref{eq: pm} is fulfilled.

\subsection{Remaining mass and entropy}

In this case we take sieves of the form 
$\BB_T = R \dot H_1^{\gamma+1/2}(\TT) + \eps_T\dot L_1^2(\TT)$,
where $\gamma$ and $R$ will be determined below. 

For the remaining mass we have
\[
\Pi(b \not\in \BB_T) \le \int_0^\gamma \lambda_T(\alpha)\,d\alpha + 
\int_\gamma^\infty \PP(W^{\alpha, 1} \not \in \BB_T)\lambda_T(\alpha)\,d\alpha.
\]
For $\alpha\ge\gamma$  we have $\BB_T \supset R \dot H_1^{\alpha+1/2}(\TT) + \eps_T\dot L_1^2(\TT)$.
Hence, by the Borell-Sudakov inequality, 
\[
\PP(W^{\alpha, 1} \not \in \BB_T) \le 1-\Phi(\Phi^{-1}(\PP(\|W^{\alpha, 1}\|_2 \le \eps_T)) + R).
\]
Note that $\|W^{\alpha,1}\|_2\le \|W^{\gamma,1}\|_2,$ so $\PP(\|W^{\alpha, 1}\|_2 \le \eps_T) \ge \PP(\|W^{\gamma, 1}\|_2 \le \eps_T)$.
By Lemma \ref{lem:smallballbalphaL}, the latter is bounded from below by $\exp(-C_\gamma\eps_T^{-1/\gamma})$
for a $C_\gamma > 0$.
We note that $C_\gamma$ depends continuously on $\gamma$, through the continuous function $f$ in Lemma 
\ref{lem:smallballbalphaL}. Below we will chose $\gamma$ to be in a shrinking neighbourhood of $\beta$, 
which is fixed. Hence, for this choice of $\gamma$, we have that 
$\PP(\|W^{\gamma, 1}\|_2 \le \eps_T) \ge \exp(-C\eps_T^{-1/\gamma})$ for a constant $C> 0$ that is 
independent of $\gamma$.
We conclude that for $\gamma\le\alpha$, 
\[
\PP(W^{\alpha, 1} \not \in \BB_T) \le \exp\left(-\left({R}-\sqrt{C'\Big(\frac{1}{\eps_T}\Big)^{1/\gamma}}\right)^2\right) 
\]
for some $C' > 0$. Taking $R$ a large multiple of $\eps_T^{-1/(2\gamma)}$ this is bounded by $\exp(-D\eps_T^{-1/\gamma})$
for a given constant $D > 0$. 
For the other term, observe that by definition of $\lambda$, 
\[
\int_0^\gamma \lambda(\alpha)\,d\alpha \le  \gamma C_T^{-1}e^{-T^{1/(1+2\gamma)}}\le \gamma e^{-T^{1/(1+2\gamma)}},
\]
since $C_T\ge\frac{\log T}{2\exp(e)}.$ 
Putting things together, we have 
\[
\Pi(b \not\in \BB_T) \le e^{-D\eps_T^{-1/\gamma}} +  \gamma e^{-T^{1/(1+2\gamma)}}.
\]
If we choose $\gamma = \beta/(1+C/\log T)$ for a large enough constant $C > 0$, then the
right-hand side is smaller than $\exp(-DT\eps^2_T)$, as desired. 

For the entropy we have, as before, 
\[
\log N(2\eps_T, \BB_T, \|\cdot\|_2) \le C \Big(\frac R{\eps_T}\Big)^{2/(1+2\gamma)}.
\]
For the choice of $R$ that we made the right side is a constant times $\eps_T^{-1/\gamma}$, 
which by the choice of $\gamma$ is bounded by a constant times $T\eps^2_T$.

\appendix

\section{Appendix}

\subsection{Unique weak solution of the periodic SDE}

\begin{lemma}\label{lemma:uniquesolutiontosde}
For $b \in\dot L^2(\T),$ the SDE  \eqref{eq: sde} has a unique weak solution.
\end{lemma}

\bigskip

\begin{proof}
Note that condition (ND) of \cite[Theorem 5.15]{karatzas} holds. Since for $0<\eps<1/2$ and $x\in\re$ we have 
$$\int_{x-\eps}^{x+\eps}|b(x)|dx\le \int_0^1|b(x)|\cdot 1dx\le \|b\|_2<\infty,$$ also condition (LI) of the theorem holds. Thus there exists a unique weak solution up to an explosion time. We will show that a solution to \eqref{eq: sde} is not explosive with probability 1. We do this by proving that the conditions of \cite[Proposition 5.22]{karatzas} are satisfied. Note that condition (ND)' holds. Furthermore for $0<\eps<1/2$ and for all $x\in\re$ we have 
$$\int_{x-\eps}^{x+\eps}(1+|b(y)|)dy\le\int_0^1(|b(y)|+1)\cdot1\le\||b|+1\|_2\le\|b\|_2+1<\infty,$$
thus condition (LI)' also holds. Define 
$$s(x)=\int_0^x\exp\left\{-2\int_0^\xi b(\zeta)d\zeta\right\}d\xi.$$ Since $b$ is 1-periodic and has mean zero, it follows that for all $x\in[0,1)$ and all $k\in\Z,s(x+k)=s(x)+k\int_0^1\exp\left\{-2\int_0^\xi b(\zeta)d\zeta\right\}d\xi,$ thus $s(x)\to\pm\infty,$ as $x\to\pm\infty.$ 
Hence \cite[Proposition 5.22]{karatzas} implies $-\infty<X_t<\infty$ almost surely, for all $t\in\re.$ 
This completes the proof.
\end{proof}

\subsection{The measures $P^{b}$  are all equivalent}

\begin{lemma}\label{lemma:equivalentmeasures}
For every $T>0$ and $b_1,b_2\in\dot L^2(\T),$ the measures $P^{b_1}=P^{b_1,T}$ and $P^{b_2}=P^{b_2,T}$ are equivalent. 	
\end{lemma}

\bigskip

\begin{proof}
Fix $T > 0$ and $b_0 \in \dot L^2(\T)$.
For every $b \in \dot L^2(\T)$ we have the occupation times formula
\[
\int_0^Tb^2(X_s)ds=\int_0^1b^2(x)L_T^\circ(x;X)\,dx.
\]
Since $P^{b_0}$-a.s.\ we have
$\|L^\circ_T(\cdot; X)\|_\infty < \infty$, it follows that 
for every $b \in \dot L^2(\T)$, we have  $\int_0^Tb^2(X_s)\,ds<\infty$, a.s.\ with respect to $P^{b_0}$.
Hence, by  Theorem III.5.38 of \cite{jacod},  all measures $P^{b, T},b\in\dot L^2(\T)$, are equivalent.
\end{proof}

\subsection{Proof of \Cref{lem:posterioriswelldefined}}\label{proof:posteriorwelldefined}

(i). 
We deal with the Lebesgue integral and the stochastic integral in \eqref{eq: lik} separately.
First note that by the occupation times formula, $\int_0^T b^2(X_t)dt = \int_0^1 b^2(x)L^\circ_T(x; X)\,dx$. Since
$P^{b_0}$-a.s.\ we have
$\|L^\circ_T(\cdot; X)\|_\infty < \infty$, this implies that $b \mapsto \int_0^T b^2(X_t)dt$
is a continuous and hence measurable functional on $\dot L^2(\TT)$. 

Using the SDE for $X,$ the stochastic integral in \eqref{eq: lik} can be written as the sum of a 
Lebesgue integral and a Brownian integral. The Lebesgue integral can be handled as in the preceding 
paragraph. To show that the Brownian integral $b \mapsto \int_0^T b(X_t)\,dB_t$ is measurable on $\dot L^2(\TT)$
we write 
\[
\dot L^2(\TT) = \bigcup_{K \in \NN} B_K,
\]
where $B_K = \{b \in \dot L^2(\TT): \|b\|_2 \le K\}$.
On every ball $B_K$ the measurability follows from 
the first statement of the Stochastic Fubini theorem as given in Theorem 2.2 of \cite{Veraar}. 
Indeed,  condition (2.1) of \cite{Veraar} translates into the requirement that, $P^{b_0}$-a.s.,  
\[
\int_{B_K} \Big(\int_0^T b^2(X_t)\,dt\Big)^{1/2}\, \Pi(db) < \infty. 
\]
This is clearly fulfilled since, by the occupation times formula again, the left-hand side is bounded by 
$K\|L^\circ_T(\cdot; X)\|_\infty^{1/2}$.

(ii). 
For the upper bound we note that the $P^{0,T}$-expectation of the {denominator} equals $1$, 
hence it is $P^{0,T}$-a.s.\ finite. But then also $P^{b_0}$-a.s., since the measures
are equivalent by Lemma \ref{lemma:equivalentmeasures}.

For the lower bound we first observe that
since $\Pi$ is  probability measure on $\dot L^2(\TT)$ there exists a $K> 0$ such that 
$\Pi(B_K) > 0$. Let $\tilde \Pi$ be the restriction of $\Pi$ to $B_K$, renormalised so that it is a
probability measure again. Then it follows from Jensen's inequality that 
\begin{align*}
\int p_b(X^T)\,\Pi(db) & \ge \Pi(B_K) \int p_b(X^T)\,\tilde\Pi(db)\\
& \ge \Pi(B_K) \exp\Big(\int \log  p_b(X^T)\,\tilde\Pi(db)\Big).
\end{align*}
Hence, it suffices to show that $P^{b_0}$-a.s., 
\[
\Big|\int \log  p_b(X^T)\,\tilde\Pi(db)\Big| < \infty.
\]
As before the log-likelihood can be written as a sum of Lebesgue and stochastic integrals.
Dealing with the Lebesgue integrals is straightforward, in view of the occupation times
formula again and the a.s.\ finiteness of $\|L_T^\circ(\cdot; X)\|_\infty$. 
It remains to show that 
$P^{b_0}$-a.s., 
\[
\Big|\int \Big(\int_0^T b(X_t)\,dW_t\Big)\,\tilde\Pi(db)\Big| < \infty.
\]
But this follows from the stochastic Fubini theorem of  \cite{Veraar} again, since as shown 
above the necessary condition for the theorem to hold is fulfilled.

\bigskip

\bibliographystyle{harry}
\bibliography{Bieb}

\end{document}